\documentclass[11pt]{amsart}
\usepackage{graphicx}
\usepackage{amsmath}
\usepackage{amsfonts}
\usepackage{amssymb}
\theoremstyle{plain}
\newtheorem{theorem}{Theorem}[section]
\newtheorem{lemma}[theorem]{Lemma}
\newtheorem{remark}[theorem]{Remark}
\numberwithin{equation}{section}

\begin{document}
\title[Linearity of homogeneous order one solutions]{Linearity of homogeneous order one solutions to elliptic equations in
dimension three}
\author[Han]{Qing Han}
\address{Department of Mathematics\\
Notre Dame University\\
Notre Dame, IN 46556}
\email{qhan@nd.edu}
\author[Nadirashvili]{Nikolai Nadirashvili}
\address{Department of Mathematics\\
University of Chicago\\
5734 S. University Ave., Chicago, IL 60637}
\email{nicholas@math.uchicago.edu}
\author[Yuan]{Yu Yuan}
\address{Department of Mathematics\\
University of Chicago\\
5734 S. University Ave., Chicago, IL 60637\\
and University of Washington\\
Seattle, WA 98195}
\email{yuan@math.uchicago.edu}
\thanks{All authors are partially supported by NSF grants, and the first author also
by Sloan Research Fellowship}
\date{October 29, 2001. Submitted on November 12, 2001.}
\maketitle

\section{\textbf{Introduction}}

In this note we prove that any homogeneous order one solution to
non-divergence elliptic equations in $\mathbb{R}^{3}$ must be linear. Consider
the general equations in $\mathbb{R}^{n}$
\begin{equation}
\sum_{i,j=1}^{n}a_{ij}(x)D_{ij}u=0, \label{ND}%
\end{equation}
where the coefficients satisfy
\[
\lambda I\leq(a_{ij}(x))\leq\lambda^{-1}I
\]
for some positive constant $\lambda>0$, and the dimension $n\geq3.$ Safonov
constructed homogeneous order $\alpha$ solutions, with $\alpha\in(0,1)$, to
(\ref{ND}) in \cite{S1}, where he showed the unimprovability of the estimates
of H\"{o}lder exponent for solutions to (\ref{ND}) by Krylov and himself. The
homogeneity $\alpha$ with $\alpha<1$ plays an essential role in the
construction. Later on, Safonov asked in [11, p.49] whether one can construct
nontrivial homogeneous order one solutions to (\ref{ND}). Our result indicates
that it is impossible to do so in $\mathbb{R}^{3}.$

\begin{theorem}
\label{main} Any homogeneous order one strong solution $u$ to (\ref{ND}) in
$\mathbb{R}^{3}$ with $u\in W_{loc}^{2,2}\left(  \mathbb{R}^{3}\right)  $ must
be a linear function.
\end{theorem}

On the other hand, one does have nontrivial homogeneous order one solutions to
(\ref{ND}) in $\mathbb{R}^{4}\backslash\{0\}.$ In fact let $\left(  x,f\left(
x\right)  \right)  \in\mathbb{R}^{n}\times$ $\mathbb{R}^{k}$ be a
nonparametric minimal surface. Then each component of $f$ satisfies
(\ref{ND}), with $\left(  a_{ij}(x)\right)  $ being the inverse of the induced
metric of the minimal surface in $\mathbb{R}^{n+k}$ (cf. [7, p.3]). Through
Hopf fibration, Lawson and Osserman constructed a minimal cone $\left(
x,f\left(  x\right)  \right)  \in\mathbb{R}^{4}\times$ $\mathbb{R}^{3},$
where
\[
f\left(  x\right)  =\frac{\sqrt{5}}{2}\frac{1}{\left|  x\right|  }\left(
x_{1}^{2}+x_{2}^{2}-x_{3}^{2}-x_{4}^{2},2x_{1}x_{3}+2x_{2}x_{4},2x_{2}%
x_{3}-2x_{1}x_{4}\right)
\]
(see [7, Theorem 7.1]). Now each component of $f$ is a desired nontrivial
solution to (\ref{ND}) in $\mathbb{R}^{4}\backslash\{0\}.$ Actually by
noticing that the graph of $u=\left(  x_{1}^{2}+x_{2}^{2}-x_{3}^{2}-x_{4}%
^{2}\right)  /\left|  x\right|  $ is a saddle surface, one can easily
construct coefficients $a_{ij}(x)$ in $\mathbb{R}^{4}\backslash\{0\}$ so that
$u$ satisfies (\ref{ND}) in $\mathbb{R}^{4}\backslash\{0\}.$

Theorem \ref{main} gives a simple ``PDE'' proof of a well-known result
obtained by many authors in 1970's, which states that any nonparametric
minimal cone of dimension three must be flat. Let $\left(  x,f\left(
x\right)  \right)  \in\mathbb{R}^{3}\times$ $\mathbb{R}^{k}$ be the minimal
cone with $f\left(  tx\right)  =tf\left(  x\right)  ,\,f$ $\in C^{\infty
}\left(  \mathbb{R}^{3}\backslash\{0\}\right)  ,$ then each component of $f$
satisfies (\ref{ND}). By Theorem 1.1, $f$ must be linear, or the minimal cone
is flat.

From Theorem 1.1 one also sees that any smooth homogeneous order two solution
in $\mathbb{R}^{3}\backslash\{0\}$ to the fully nonlinear elliptic equation
$F\left(  D^{2}u\right)  =0$ must be a quadratic polynomial. To get this
conclusion, one simply apply Theorem 1.1 to the gradient $\nabla u$ and
(\ref{ND}) with $a_{ij}(x)=\frac{\partial F}{\partial m_{ij}}\left(
D^{2}u\right)  .$ In contrast, the second author constructed a non-quadratic
homogeneous order two solution to some equation $F\left(  D^{2}u\right)  =0$
in $\mathbb{R}^{12},$ which provides a counterexample to the regularity for
fully nonlinear elliptic equations (see \cite{N}).

The heuristic idea of the proof of Theorem 1.1 is the following simple
geometric observation. The closed saddle surface $\nabla u\left(
\mathbb{S}^{2}\right)  $ must be a point. Therefore, $u$ is linear. More
precisely, one consider the surface $\Sigma$ parametrized by the gradient
$\nabla u:\mathbb{S}^{2}\rightarrow\mathbb{R}^{3}.$ Because of (\ref{ND}),
$\Sigma$ is a saddle surface at $\nabla u\left(  x\right)  $ with
$D^{2}u\left(  x\right)  \neq0.$ However the supporting plane with normal $x$
touches $\Sigma$ at $\nabla u\left(  x\right)  .$ Thus $D^{2}u\left(
x\right)  \equiv0$ and it follows that $u$ is linear, see Section 2.

The ideas of gradient map and supporting plane are already in an early paper
by Alexandrov \cite{A}. In fact, under the assumption that the homogeneous
order one function $u$ is analytic in $\mathbb{R}^{3}\backslash\{0\}$ and the
Hessian $D^{2}u$ is either non-definite or $0$ at each point, Alexandrov showed
that $u$ must be a linear function. Roughly speaking, if $u$ is analytic in
$\mathbb{R}^{3}\backslash\{0\}$, the set $\mathcal{S}(u)=\left\{
x\in\mathbb{S}^{2}|D^{2}u\left(  x\right)  =0\right\}  $ is either isolated or
the whole $\mathbb{S}^{2}$. Alexandrov proved that the supporting plane to
$\Sigma$ is unique at $\nabla u\left(  x\right)  $ with $x$ being an isolated
point of $\mathcal{S}(u).$ Since the surface $\Sigma$ has supporting planes
with normal along all the directions in $R^{3},$ Alexandrov excluded case of
$\mathcal{S}(u)$ being isolated and proved the result. It is interesting to
note that the concept of gradient maps and supporting planes was further
employed in the later development of the Alexandrov-Bakel'man-Pucci maximum principle.

When the coefficients $a_{ij}$ are $C^{\alpha}$ and the solution $u$ is
$C^{2,\alpha}$ away from the origin, we can show that the set $\mathcal{S}(u)$
is either isolated or the whole $\mathbb{S}^{2}.$ Coupled with Alexandrov's
argument, one sees that Theorem 1.1 holds in the $C^{2,\alpha}$ setting. After
the work in $C^{2,\alpha}$ case was done, we found that Pogorelov \cite{P}
generalized the above Alexandrov's result to $C^{2}$ functions. The argument
in \cite{P} is more involved. Since our approach in $C^{2,\alpha}$ case is 
interesting in
its own right and short, we also include it here in Section 3.

The authors wish to thank L. Nirenberg for informing and mailing them the
inspiring reference \cite{A} .

\section{\textbf{Proof of the Main Theorem}}

Let $h\left(  x_{1},x_{2}\right)  =u\left(  x_{1},x_{2},1\right).$  Then the
homogeneous order one function $u\left(  x_{1},x_{2},x_{3}\right)
=x_{3}h\left(  \frac{x_{1}}{x_{3}},\frac{x_{2}}{x_{3}}\right)  ,$ and the gradient
$\nabla u$ and the Hessian $D^{2}u$ have the following representation%
\begin{equation}
\nabla u\left(  x_{1},x_{2},1\right)  =\left(  h_{1},h_{2},h-x_{1}h_{1}%
-x_{2}h_{2}\right)  , \label{repg}%
\end{equation}%
\begin{equation}
D^{2}u\left(  x_{1},x_{2},1\right)  =\left[
\begin{array}
[c]{ccc}%
1 & 0 & 0\\
0 & 1 & 0\\
-x_{1} & -x_{2} & 1
\end{array}
\right]  \left[
\begin{array}
[c]{ccc}%
h_{11} & h_{12} & 0\\
h_{21} & h_{22} & 0\\
0 & 0 & 0
\end{array}
\right]  \left[
\begin{array}
[c]{ccc}%
1 & 0 & -x_{1}\\
0 & 1 & -x_{2}\\
0 & 0 & 1
\end{array}
\right]  . \label{repH}%
\end{equation}
From our assumption that the homogeneous order one solution $u$  is in $W_{loc}%
^{2,2}\left(  \mathbb{R}^{3}\right)  ,$ it follows that $h\in W_{loc}%
^{2,2}\left(  \mathbb{R}^{2}\right)  $ is a strong solution to%
\begin{equation}
\sum_{i,j=1}^{2}A_{ij}\left(  x_{1},x_{2}\right)  D_{ij}h=0, \label{ND2d}%
\end{equation}
where the coefficients $A_{ij}$ are in terms of $a_{ij}\left(  x_{1}%
,x_{2},1\right)  ,\;x_{1},\;x_{2}$ and satisfy the ellipticity condition with
some $\lambda\left(  x_{1},x_{2}\right)  .$

\begin{lemma}
\label{support}For any $\nu\in\mathbb{S}^{2}$, the supporting plane with
normal $\nu$ must touch the surface $\Sigma$ parametrized by $\nabla
u:\mathbb{S}^{2}\rightarrow\mathbb{R}^{3}$ at $\nabla u\left(  \nu\right)  $
or $\nabla u\left(  -\nu\right)  .$
\end{lemma}

\begin{proof}
First we notice that $\nabla u\in C^{\alpha}\left(  \mathbb{S}^{2}\right)  .$
In fact each component of $\left(  h_{1},h_{2}\right)  ,$ say $h_{1}$ is a
$W_{loc}^{1,2}\left(  \mathbb{R}^{2}\right)  $ weak solution to a divergence
equation%
\[
\sum_{i,j=1}^{2}D_{i}\left(  B_{ij}D_{j}h\right)  =0,
\]
where $B_{11}=A_{11}/A_{22},\;B_{12}=2A_{12}/A_{22},\;B_{22}=1.$ It follows
that $h\in C^{1,\alpha}\left(  S^{2}\right)  $ (cf. [5, p,284-285 or Theorem
12.4]). By (\ref{repg}) we see $\nabla u\in C^{\alpha}\left(  \mathbb{S}%
^{2}\right)  $ with $\alpha$ depending only on the original ellipticity
$\lambda.$ Therefore, for any $\nu\in\mathbb{S}^{2}$, the supporting plane
$P_{\nu}$ with normal $\nu$ must touch the surface $\Sigma$ somewhere
($\Sigma$ is on the opposite side of $\nu$).

Without loss of generality,we assume $\nu=\left(  1,0,0\right)  .$

\textbf{Claim:} $P_{\left(  1,0,0\right)  }$ touches $\Sigma$ at $\nabla
u\left(  x_{1},x_{2},0\right)  $ with $\left(  x_{1},x_{2},0\right)
\in\mathbb{S}^{2}.$

Suppose $P_{\left(  1,0,0\right)  }$ touches $\Sigma$ at $\nabla u\left(
x_{1},x_{2},x_{3}\right)  $ with $x_{3}\neq0,$ say $x_{3}>0.$ Then
$u_{1}\left(  x_{1},x_{2},1\right)  =u_{1}\left(  \frac{x_{1}}{x_{3}}%
,\frac{x_{2}}{x_{3}},1\right)  =h_{1}\left(  \frac{x_{1}}{x_{3}},\frac{x_{2}%
}{x_{3}}\right)  $ would achieve its maximum at $\left(  \frac{x_{1}}{x_{3}%
},\frac{x_{2}}{x_{3}}\right)  .$ Because of (\ref{ND2d}), $h_{1}$ satisfies
the strong maximum principle (cf. [5,Theorem 8.19]). We then have a
contradiction unless $h_{1}\equiv const..$ In the latter case, $u_{1}\left(
x_{1},x_{2},x_{3}\right)  \equiv const.$ for $x_{3}>0,$ our claim holds.

Similarly, by applying the above argument to $u\left(  x_{1},1,x_{3}\right)
,$ we see that $P_{\left(  1,0,0\right)  }$ also touches $\Sigma$ at $\nabla
u\left(  x_{1},0,x_{3}\right)  $ with $\left(  x_{1},0,x_{3}\right)
\in\mathbb{S}^{2}.$ Therefore, $P_{\left(  1,0,0\right)  }$ must touch
$\Sigma$ at $\nabla u\left(  1,0,0\right)  $ or $\nabla u\left(
-1,0,0\right)  .$
\end{proof}

We now present the proof of Theorem \ref{main}.

\begin{proof}
By our assumption on $u,$ $u\in W^{2,2}\left(  \mathbb{S}^{2}\right)  .$ If
$D^{2}u=0$ almost everywhere on $\mathbb{S}^{2},$ then $u$ is already linear.
Otherwise, we pick a Lebesgue point $x^{\ast}\in\mathbb{S}^{2}$ for $D^{2}u$
with $D^{2}u\left(  x^{\ast}\right)  \neq0,$ say $x^{\ast}=\left(
1,0,0\right)  .$ We may also assume $x^{\ast}$ is a Lebesgue point for $\nabla
u$ and $\left(  a_{ij}\left(  x\right)  \right)  .$ By Lemma \ref{support},
the supporting planes $P_{\left(  1,0,0\right)  }$ and $P_{\left(
-1,0,0\right)  }$ touches $\Sigma$ at $\nabla u\left(  1,0,0\right)  $ or
$\nabla u\left(  -1,0,0\right)  .$ If both planes touch $\Sigma$ at the same
point, then we see that $u_{3}\equiv c.$ Consequently the homogeneous order
one function $v\left(  x_{1},x_{2}\right)  =u-cx_{3}$ satisfies (\ref{ND}) and
it follows that $v$ is linear, or $u$ is also linear.

So we are left with the case that $P_{\left(  1,0,0\right)  }$ and $P_{\left(
-1,0,0\right)  }$ touch $\Sigma$ at different point. We may assume $P_{\left(
1,0,0\right)  }$ touches $\Sigma$ at $\nabla u\left(  1,0,0\right)  .$ It
means that%
\begin{equation}
u_{3}\left(  x_{1},x_{2},1\right)  =h-x_{1}h_{1}-x_{2}h_{2}\leq h\left(
0\right)  \;\;\;\text{near\ \ }\left(  0,0\right)  . \label{Sbp}%
\end{equation}
By (\ref{repg}) and (\ref{repH}), $\left(  0,0\right)  $ is a Lebesgue point
for $D^{2}h$ and $\nabla h.$ It follows that (cf. [4, Appendix C])%
\[
h\left(  x_{1},x_{2}\right)  =h\left(  0\right)  +h_{1}\left(  0\right)
x_{1}+h_{2}\left(  0\right)  x_{2}+\frac{1}{2}\sum_{i,j=1}^{2}h_{ij}\left(
0\right)  x_{i}x_{j}+0\left(  \left|  x\right|  ^{2}\right)  .
\]
Because of (\ref{ND2d}) and $D^{2}h\left(  0,0\right)  \neq0,$ we may assume
$h_{11}\left(  0\right)  =a,\;h_{12}\left(  0\right)  =0,\;h_{22}\left(
0\right)  =-b$ with $a,b>0.$ A simple computation yields $u\left(
0,x_{2},1\right)  =h\left(  0\right)  +\frac{1}{2}bx_{2}^{2}+o\left(
x_{2}^{2}\right)  $ which contradicts (\ref{Sbp}). This completes the proof.
\end{proof}

\section{\textbf{Another Proof in the $C^{2,\alpha}$ Case}}

In this section, we present yet another proof of Theorem \ref{main} in the
$C^{2,\alpha}$ case. Suppose $u:\mathbb{R}^{3}\setminus\{0\}\rightarrow
\mathbb{R}$ is a $C^{2,\alpha}$ homogeneous order one function $u=rg\left(
\theta_{1},\theta_{2}\right)  ,$ where $\left(  r,\theta_{1},\theta
_{2}\right)  $ is the spherical coordinates with $x_{1}=r\cos\theta_{2}%
\cos\theta_{1},$ $x_{2}=r\cos\theta_{2}\sin\theta_{1},$ $x_{3}=r\sin\theta
_{2}.$ Then the Hessian $D^{2}u=\left(  \frac{\partial^{2}u}{\partial
x_{i}\partial x_{j}}\right)  $ has the following representation%
\begin{equation}
D^{2}u\left(  x\right)  =\frac{1}{r}\mathbf{R}\left(  \theta_{1},\theta
_{2}\right)  H\mathbf{R}^{^{\prime}}\left(  \theta_{1},\theta_{2}\right)
\label{Hs}%
\end{equation}
where
\[
H=\left[
\begin{array}
[c]{ccc}%
0 & 0 & 0\\
0 &\frac{1}{\cos^{2}\theta_{2}} 
\frac{\partial^{2}g}{\partial\theta_{1}^{2}}-\tan\theta_{2}\frac{\partial
g}{\partial\theta_{2}}+g & \frac{\partial^{2}g}{\partial\theta_{1}%
\partial\theta_{2}}+\frac{\sin\theta_{2}}{\cos^{2}\theta_{2}}\frac{\partial
g}{\partial\theta_{1}}\\
0 & \ast & \frac{\partial^{2}g}{\partial\theta_{2}^{2}}+g
\end{array}
\right]
\]
and $\mathbf{R}\left(  \theta_{1},\theta_{2}\right)  $ is the rotation from
$\left(  \partial x_{1,}\partial x_{2},\partial x_{3}\right)  $ to $\left(
\partial r,\frac{1}{r}\partial\theta_{1},\frac{1}{r}\partial\theta_{2}\right)
.$ We see that $D^{2}u=0$ if and only if $H=0.$ We define the singular set
$\textsl{S}\left(  u\right)  $ as
\[
\textsl{S}(u)=\{x\in\mathbb{S}^{2}|D^{2}u(x)=0\}.
\]

\begin{lemma}
\label{sing} Let $u$ be as in Theorem \ref{main}. Then one of the following holds:

(1) $\textsl{S}(u)$ is empty;

(2) $\textsl{S}(u)$ consists of finitely many points;

(3) $\textsl{S}(u)=$ $\mathbb{S}^{2}.$
\end{lemma}

\begin{proof}
Suppose $\textsl{S}(u)$ is not empty. For any $p\in\textsl{S}(u),$ we prove that
either $p$ is isolated or $\textsl{S}(u)$ contains a neighborhood of $p$ on
$\mathbb{S}^{2}.$ Without loss of generality, we assume $p=\left(
1,0,0\right)  ,$ $a_{ij}\left(  p\right)  =\delta_{ij}.$ Note that
$\textsl{S}(u)=\textsl{S}(u+ax_{1}+bx_{2}+cx_{3}),$ we may also assume $\nabla
u\left(  p\right)  =0,$ then $u\left(  p\right)  =p\cdot\nabla u\left(
p\right)  =0.$ Correspondingly, we have
\[
\left(  \frac{\partial^{2}g}{\partial\theta_{i}\partial\theta_{j}}\right)
\left(  0,0\right)  =0,\;\;\left(  \frac{\partial g}{\partial\theta_{1}},
\frac{\partial g}{\partial\theta_{2}}\right)  \left(  0,0\right)
=0,\;\;\text{and}\;\;g\left(  0,0\right)  =0.
\]
Also $g$ satisfies
\[
\sum_{i,j=1}^{2}A_{ij}\left(  \theta_{1},\theta_{2}\right)  \frac{\partial
^{2}g}{\partial\theta_{i}\partial\theta_{j}}+\sum_{i=1}^{2}B_{i}\left(
\theta_{1},\theta_{2}\right)  \frac{\partial g}{\partial\theta_{i}}+C\left(
\theta_{1},\theta_{2}\right)  g=0,
\]
where the $C^{\alpha}$ coefficients $\left(  A_{ij}\right)  $ satisfy the
ellipticity condition with the same constant $\lambda,$ $A_{ij}\left(
0,0\right)  =\delta_{ij},$ $B_{i}$ and $C$ are $C^{\alpha}$ functions near
$\left(0,0\right)$.
Suppose $g$ vanishes up to infinity order at $\left(  0,0\right)  ,$ then by
the Carleman unique continuation (cf. [3, p.124]), $g\equiv0$ in a
neighborhood of
$\left(  0,0\right)  .$ It follows that $D^{2}u=0$ on $\mathbb{S}^{2},$ since
$u$ is $C^{2}$.
Suppose $g$ vanishes up to order $k-1$ at $\left(  0,0\right)  .$ By our
assumption $k\geq3.$ We apply the result of [2, Theorem 1] to obtain
\[
g=P\left(  \theta_{1},\theta_{2}\right)  +R\left(  \theta_{1},\theta
_{2}\right)  ,
\]
where the homogeneous order $k$ polynomial $P$ satisfies $\sum_{i,j=1}^{2}
A_{ij}\left(  0,0\right)  \frac{\partial^{2}P}{\partial\theta_{i}
\partial\theta_{j}}=0,$ or $\triangle P=0,$ and the remainder satisfies
\[
R\sim O\left(  \left|  \theta\right|  ^{k+\alpha}\right)  ,\;\;\nabla R\sim
O\left(  \left|  \theta\right|  ^{k-1+\alpha}\right)  ,\;\;D^{2}R\sim O\left(
\left|  \theta\right|  ^{k-2+\alpha}\right)
\]
with $\alpha\in\left(  0,1\right)  .$ It is a simple fact that $\left\{
\left(  \theta_{1},\theta_{2}\right)  |D^{2}P=0\right\}  =\left\{\left(
0,0\right)\right\}$ and consequently $\left\{  \left(
\theta_{1},\theta_{2}\right)
|D^{2}H=0\right\}  =\left\{\left(  0,0\right)\right\}  .$ Hence $p$ is an
isolated zero point
of $D^{2}u.$ Therefore $\textsl{S}(u)$ consists of finitely many points in
this case.
\end{proof}

We consider the surface $\Sigma$ parametrized by $\nabla u:\mathbb{S}
^{2}\rightarrow\mathbb{R}^{3}.$ From (\ref{Hs}), it follows that the Hessian
always has one zero eigenvalue. Let $\lambda_{1}\left(  x\right)  $ and
$\lambda_{2}\left(  x\right)  $ be the other two eigenvalues of $D^{2}u\left(
x\right)  .$ Because of equation (\ref{ND}), $\lambda_{1}\left(  x\right)
\lambda_{2}\left(  x\right)  <0$ for all $x\in\mathbb{S}^{2}\backslash
\textsl{S}(u).$

\begin{lemma}
\label{curv} For any $x\in\mathbb{S}^{2}\backslash\textsl{S}(u)$, the surface
$\Sigma$ is $C^{2,\alpha}$ at $\nabla u(x)$ with a normal vector given by $x$
and the two principle curvatures given by $-1/\lambda_{1}\left(  x\right)  $
and $-1/\lambda_{2}\left(  x\right)  .$
\end{lemma}

\begin{proof}
We may assume $x=p=(0,0,1)\in\mathbb{S}^{2}\setminus\textsl{S}(u).$ Then
locally at $\nabla u(p),$ $\Sigma$ can be represented by
\begin{equation}
\mathbf{F}(x_{1},x_{2})=\nabla u(x_{1},x_{2},\sqrt{1-x_{1}^{2}-x_{2}^{2}}).
\label{surf}
\end{equation}
By differentiating the identity $x\cdot\nabla u(x)=u\left(  x\right)  $ twice
with respect to $x_{i},x_{j}$ for $i,j=1,2,$ we obtain
\begin{equation}
u_{3i}(p)=0,\ u_{3ij}(p)=-u_{ij}(p).\label{value}
\end{equation}
By differentiating (\ref{surf}) and making use of (\ref{value}), we have
\[
\mathbf{F}_{i}(0)=(u_{i1}(p),u_{i2}(p),0)
\]
and
\[
\mathbf{F}_{ij}^{(3)}(0)=-u_{ij}(p)
\]
where $\mathbf{F}_{ij}^{(3)}$ denotes the third component of the vector
$\mathbf{F}_{ij}$. Then we get
\[
\mathbf{F}_{1}(0)\times\mathbf{F}_{2}(0)=(0,0,(u_{11}u_{22}-u_{12}^{2})(p)).
\]
It is easy to see that $(u_{11}u_{22}-u_{12}^{2})(p)=\lambda_{1}\lambda
_{2}(p)\neq0.$ Hence $(0,0,1)$ is a normal vector for $\Sigma$ at $\nabla
u(p).$ We see that $\nabla u =G^{-1}$ near $(0,0,1),$ where $G^{-1}$ is the
inverse of the Gauss map of $\Sigma.$
Thus  $\Sigma$ is $C^{2,\alpha}$ nearby. We also get the first and second
fundamental forms of $\Sigma$ at $\nabla u(p)$ as follows
\begin{align*}
I &  =(u_{11}^{2}+u_{12}^{2})dx_{1}^{2}+2u_{12}(u_{11}+u_{22})dx_{1}
dx_{2}+(u_{21}^{2}+u_{22}^{2})dx_{2}^{2},\\
II &  =-u_{11}dx_{1}^{2}-2u_{12}dx_{1}dx_{2}-u_{22}dx_{2}^{2},
\end{align*}
where all $u_{ij}$ are evaluated at $p$. Therefore the two principle
curvatures are $-1/\lambda_{1}$ and $-1/\lambda_{2}.$
\end{proof}

\begin{remark}
In computing the principle curvatures of  $\Sigma$ in terms of the eigenvalues
of $D^{2}u,$ we differentiate $u$ three times. By approximation, the
conclusion still holds for $u\in C^{2,\alpha}.$
\end{remark}

Now we are ready for another proof of Theorem {\ref{main} in the }%
$C^{2,\alpha}$ case.

\begin{proof}
We prove that $\textsl{S}(u)=\mathbb{S}^{2}$ by excluding the case (1) and (2)
in Lemma \ref{sing}. We assume $\textsl{S}\left(  u\right)  $ consists of at
most finitely many points. First, the gradient surface $\Sigma$ has supporting
planes with normals along all the directions in $\mathbb{R}^{3}$. However, the
saddle points in $\left\{  \nabla u\left(  x\right)  |D^{2}u\left(  x\right)
\neq0\right\}  $ cannot support any supporting planes. Hence there are at most
finitely many points on $\Sigma$ with the supporting planes.

\noindent \textbf{Claim.}(Alexandrov) For any $x_{0}\in\textsl{S}(u)$, the 
supporting plane to
$\Sigma$ at $\nabla u(x_{0})$ can only have normal direction $x_{0}.$

Therefore there are at most finitely many supporting planes to $\Sigma$. This
is a contradiction.
Now we prove the claim. Suppose there is another supporting plane $P_{1}$ to
$\Sigma$ at $\nabla u(x_{0})$ with normal direction $x_{1}\neq x_{0}.$ Take a
small neighborhood $U$ of $x_{0}$ on $\mathbb{S}^{2}$ such that $D^{2}u\left(
x\right)  \neq0$ for any $x\neq x_{0}\in U$ and $U\cap S^{\ast}=\phi$, with
$S^{\ast}$ being a great circle through $x_{1}$ and $-x_{1}.$ This can be done
since $x_{0}$ is an isolated point in $\mathcal{S}(u)$. Now lift $P_{1}$ along
the $x_{1}$ direction to $P$ so that $P\cap\nabla u\left(  U\right)  =C$ is a
smooth close curve on $P.$ We can take a point on $C,$ say $\nabla u\left(
x^{\ast}\right)  $ with $x^{\ast}\in U$ and $x^{\ast}\neq x_{0}$, such that
the normal to the plane curve $C$ at $\nabla u\left(  x^{\ast}\right)  $ is
along the intersection $P$ and another plane through $S^{\ast}.$ Then we see
that the normal of the surface $\Sigma$ at the regular point $\nabla u\left(
x^{\ast}\right)  $ must be on $S^{\ast}.$ On the other hand, the normal at
$\nabla u\left(  x^{\ast}\right)  $ is $x^{\ast}\in U.$ This contradiction
completes the proof of the Claim.

By excluding the case (1) and (2) in Lemma
\ref{sing}, we are left with the case (3) $\mathcal{S}(u)=\mathbb{S}^{2}.$
That is $D^{2}u\equiv0$, and hence $u$ is linear. This finishes the proof of
Theorem \ref{main} in the $C^{2,\alpha}$ case.
\end{proof}


\begin{thebibliography}{99}
\bibitem{A}Alexandrov, A. D., \emph{Sur les th\'{e}r\`{e}mes d'unicit\'{e}
pour les surfaces ferme\'{e}s,} C.\ R. Acad. Sci. URSS, N.s. \textbf{22}
(1939), 99--102.

\bibitem {Bers}Bers, L., \emph{Local behavior of solution of general linear
elliptic equations}, Comm. Pure Appl. Math., \textbf{8} (1955), 473--496.

\bibitem {BL}Bers, L. and Nirenberg, L, \emph{On a representation theorem for
linear elliptic systems with discontinuous coefficients and its application},
Convegno Internazionale sulle Equatzioni Derivate e Parziali, August 1954, 111--140.

\bibitem {CCKS}Caffarelli, L. A., Crandall, M. G., Kocan, M., and \'{S}wiech,
A., \emph{On viscosity solutions of fully nonlinear equations with measurable
ingredients}, Comm. Pure Appl. Math., \textbf{49} (1996), 365--397.

\bibitem {GT}Gilbarg, D. and Trudinger, N. S., \emph{Elliptic Partial
Differential Equations of Second Order}, 2ed edition, Springer-Verlag, New
York, 1983.

\bibitem {Ha}Han, Q., \emph{Schauder estimates for elliptic operators with
applications to nodal sets}, J. Geom. Analysis, \textbf{10} (2000), 455--480.

\bibitem {LO}Lawson, H. B. Jr. and Osserman, R., \emph{Non-existence,
non-uniqueness and irregularity of solutions to the minimal surfaces system,}
Acta Math., \textbf{139} (1977), 1--17.

\bibitem {N}Nadirashvili, N., \emph{Nonclassical solutions to fully nonlinear
elliptic equations}, to appear.

\bibitem {P}Pogorelov, A. V., \emph{Solution\ of\ a problem of A. D.
Aleksandrov}, (Russian), Dokl. Akad. Nauk, \textbf{360} (1998), 317--319.

\bibitem {S1}Safonov, M. V., \emph{Unimprovabilty of estimates of Holder
constants for solutions of linear elliptic equations with measurable
coefficients,} Mat. Sb., \textbf{132} (1987), 272--288; English translation in
Math. USSR Sb., \textbf{60} (1988), 269--281.

\bibitem {S2}Safonov, M. V., \emph{Nonlinear Elliptic Equations of the Second
Order,} Lecture notes, Univ. di Firenze, 1991.
\end{thebibliography}
\end{document}